\documentclass[11pt,a4paper]{article}

\usepackage{graphicx}
\usepackage{amsmath}
\usepackage{amsfonts}
\usepackage{amssymb}
\usepackage{mathtools}
\usepackage{mathrsfs}
\usepackage{enumerate}
\usepackage{lscape}
\usepackage{longtable}
\usepackage{rotating}
\usepackage{multirow}
\usepackage{color}
\usepackage{xcolor}
\usepackage{url}
\usepackage{subfigure}
\usepackage{rotating}
\usepackage{hyperref}
\hypersetup{
    colorlinks=true, 
    linktoc=all,     
    linkcolor=blue,  
}
\usepackage{algorithmic}

\usepackage[ruled,vlined,noline,linesnumbered]{algorithm2e}

\usepackage{dsfont}

\usepackage{authblk}

\usepackage[font=scriptsize]{caption}

\usepackage{algorithmic}

\DeclareMathAlphabet \mathbfcal{OMS}{cmsy}{b}{n}

\title{Quasi-Ergodic Control of Multi-Periodic Autoregressive Processes: Formulation and Examples}
\author[1]{Vyacheslav Kungurtsev}

\affil[1]{Czech Technical University in Prague, the Czech Republic}

\begin{document}
\maketitle

\begin{abstract}
This work considers state dynamics driven by Periodic Autoregressive Moving Average noise, and control of the system over time. Such processes appear frequently in applications involving the environment, such as energy and agriculture. Managing these systems applying forecasts to make decisions that exhibit foresight and  risk aversion while maximizing profits is a challenging control problem that can be computationally difficult for standard scenario-based methods. This paper presents a formulation that explicitly enforces time-periodicity of distribution of the state, facilitating the use of periodic stochastic process basis elements as a discretization. By enforcing periodicity explicitly, an ansatz for the solution can be formed that does not require exponential scaling with time. We provide a few examples that can be modeled with this new control formulation. 
\end{abstract}

\section{Introduction}

There are a number of processes with time-dynamics that follow a stochastic process that exhibits periodicity of multiple periods. This is easily seen with anything depending on weather - there is a natural geological cycle corresponding to the sun's varying influence on the surface of the Earth across the 24 hours in a day as well as seasonal influence throughout the period of a year. Beyond this, many systems have endogenous mechanisms that can create additional periodicity through systems relying on corrective dynamics that consistently correct for overshoot and undershoot. This was first presented as a pioneering comprehensive analysis of business cycles by the economist Joseph A. Schumpeter in~\cite{schumpeter1939business}. In this work, a combination of seasonal, regular boom-bust business cycles, and long run Kondratieff technological innovation-driven cycles are overlayed to describe the time-dynamics of aggregate economic output, investment and other macroeconomic quantities.

Modeling real phenomena of such systems requires accurate and precise system identification. Naturally, there are a number of stochastic inputs into the dynamic process. This is epistemically represented through uncertainty-quantified (even if just by standard error computations) forecasts of the weather or economic spirits or other quantities unknown in the future. In statistical time series efforts to fit data to a model for forecasting, periodicity plays an important role in many such models. It appears in the modern day for, e.g., water resource management using statistical tools of Periodic Autoregressive Regression (PAR)~\cite{vecchia1985periodic} and later Periodic Autoregressive Moving Average Regression (PAMAR)~\cite{mcleod1994diagnostic}. 

The management and control of systems of this form is a challenging computational endeavor. Technically, this requires an accurate model representation of the system together with a scheme that approximates an ideal Dynamic Programming (DP) solution. In DP (see, e.g.~\cite{denardo2012dynamic}) one performs a backwards recursion from some final time, obtaining a stochastic process defined on the uncertainty of all of the realizations of the random dynamic process up to that point and parametrized by the previous-time controls. One then inputs this solution into the next-to-last time step, and continues through what is called ``backwards induction'' up until reaching the initial time instant. In practice, however, most problems are intractable as far as exact solutions.

Multistage Stochastic Programming approaches together with scenario generation, that samples future potential paths of stochastic realizations, provides a computational approximation. However, these direct discretizations suffer from a ``curse of dimensionality'' wherein exponential explosion with respect to decision and state variable dimension and stochastic discretization (that is, a finite sample of the noise at each time step) renders computational performance riddled with limitations. Yet the accuracy quickly degrades with a meaningful number of scenarios to properly quantify the uncertainty, yielding difficult choices of tradeoffs as far as the choice of model and forecast detail as well as solution algorithm to employ.

Model Predictive Control (e.g.~\cite{holkar2010overview}) is a framework for real-time decision making in physical processes. An optimal control problem defining the operation of some system for a finite horizon of future time is discretized, meaning the differential equations defining the dynamics are computed numerically at discrete time and space points. The resulting finite dimensional optimization problem is solved using appropriate software, and the optimal control choice at the first time step is implemented into the system. Subsequently the process evolves, and its state and any relevant parameters are measured or estimated, and a new optimal control problem is defined, now initialized at the current time step, with the horizon of the same length and thus extending one time step. This process is repeated, with the intention of resulting in a long term well performing and stable trajectory. 

In practical real time, one maintains an optimization problem describing the process up through a certain maximal time interval. Then, at regular intervals, one obtains novel measurements, truncates the (essentially infinite) horizon to a fixed, finite horizon, discretizes the continuous dynamics of the model, and solves a corresponidng mixed-integer non-linear optimization problem (MINLP), often simplifying the problem to form a mixed-integer linear optimization program (MILP). In control theory, this approach is known as constrained model predictive control (CMPC) of hybrid systems.

When quantities in the system are uncertain, then the optimization problem is not well defined, as the solution depends on the realization of the uncertainty. Thus one has to compute some statistical agglomeration of the objective to defined a proper optimization problem. A standard choice is taking the expectation, over the distribution of the uncertainty, and so the problem becomes to minimize the expected cost. However, such an objective is risk neutral, and an operator may be more averse to large losses than favorable to large gains. One can minimize the volatility, or standard deviation, however this would unnecessarily penalize downward volatility. A number of \emph{risk measures}~\cite{ruszczynski2006optimization} were developed in order to quantiy one-sided risk. They each satisfy certain axiomatic properties of risk measures, and different measures have their own specific features. A popular one for optimal control applications of the Conditional Value at Risk (CVaR), see, e.g.~\cite{kouri2016risk}. CVaR considers the probability mass expected beyond a certain quintile, and is parametrized by the quintle value. Thus it presents a natural interpretation of tail behavior.



In this work we study the problem of dynamic control of such processes. To this effect we consider that, at the high level, one can obtain a quasi-ergodic solution offline and then formulate transient control problems with shorter time scales in real time to ensure close-to-ergodic operation. This ensures the simultaneous long run reliability while permitting flexibility towards addressing real-time risks and opportunities. We define the problem of interest as establishing a policy such that, with respect to a given forecast of a PAR or PAMAR phenomenon, the process follows a stationary, around a moving average, multi-periodic distributional orbit. The periods will, of course, correspond to those of the statistical quantities, and the decisions will be chosen so as to realize these probabilistic constraints while minimizing a cost that includes both an expectation and the upside risk of an economic objective, usually negative profit. 

\section{Background}
Here we review the primary computational methodologies that are related to the problem that we introduce. We provide a brief description, indicate the correspondence to the problem below, and present some references.

\subsection{Dynamic Programming}
Dynamic Programming involves solving a stochastic control problem for control $u(t)$ governing the dynamics of a state $x(t)$ through the dynamics $f$ in order to minimize some operational criteria. Formally,
\[
\begin{array}{rl}
\arg\min\limits_{u(t),x(t),\,t\in[0,T]} & J(x,u):= B(x(T))+\int_0^T V(x(t),u(t))dt\\
\text{s.t. }& x\dot{x}(t) = f(x(t),u(t)),\,x(t)\in \mathcal{X},\,u(t)\in\mathcal{U}
\end{array}
\]
wherein $\mathcal{X}$ and $\mathcal{U}$ are constraint sets. This problem, to solve exactly, would require a backwards induction of $u(t)$ through the Bellman equation~\cite{bertsekas2022abstract}. In the stochastic optimal control setting, when $f$, $B$ and $V$ depend on random variables, the problem becomes considerably more challenging due to the stochastic dependence across time~\cite{bertsekas2012dynamic}. For instance $u(t_b)$ will depend on the stochastic noise that has been realized in the $\sigma$-algebra of $t_b$, and in turn influences the total objective, which influences the decision at $u(t_a)$ for $t_a<t_b$.

A popular approximation technique is Reinforcement Learning, wherein a neural network is used to define a policy of $u(t)$ given $x(t)$, while learning about dynamics~\cite{bertsekas2019reinforcement}. 

Dynamic Programming is a general formulation that encompasses many decision problems, including the one introduced in this paper. Here, we do consider a formulation with a policy solution, however, present a formulation and suggestions towards explicit computation with exact models, rather than use learning approaches.
\subsection{Multistage Stochastic Model Predictive Control}
A popular method that can be considered as a lookahead framework for Approximate Dynamic Prorgramming~\cite{powell2007approximate} is Multistage Stochastic Model Predictive Control.

Model Predictive Control (MPC) (e.g.~\cite{rawlings2012fundamentals}) performs real-time optimization by solving the optimal control up to some time horizon $H\ll T$, then applying the first control defined by the solution. The system is allowed to evolve, the true state is measured and is set as the new initial state of the next optimization problem. Through proofs of closed loop stability, MPC can be shown to facilitate a control sequence that is a good approximate solution to the long term dynamic programming.

With stochasticity, the problem becomes a Multistage Stochastic Programming Problem, which is then typically discretized stochastically with Sample Average Approximation using scenario trees. The combinatorial explosion of past-dependent scenarios is a significant challenge of this class of problems. See~\cite{mesbah2016stochastic} for a presentation and a discussion.

For the problem defined below, we intend that the solution of an optimization problem, possibly from sample average scenarios, will be the primary stream of computation. However, rather than a receding horizon approach, the approach leverages the periodicity of the noise to define a one-shot solution operation.

\subsection{Ergodic Optimization}
Ergodic optimization considers stochastic processes long run asymptotic average. Of particular importance are whether the dynamics are ergodic, that is, asymptotically converge in distribution to some stationary distribution. Generically, ergodic optimization with criteria $A(x)$ is any map $\mathcal{T}$ from a probability density to another density 
\[
\arg\min\limits_{\rho:\,\rho=\mathcal{T}}\,\int A(x) d\rho(x)
\]

In the control context, ergodicity could depend on some subset of the control policies. Formally, with $\Xi=\{\omega(t)\}$ the random noise and $\upsilon$ a randomization of policy,
\[
\arg\min\limits_{u(x,\upsilon),x(t,\xi)} \, \lim\sup\limits_{t\to\infty} \mathbb{E}_{\Xi}\left[\int_0^t V(x(\tau,\omega(0:\tau)),u(x,\upsilon)) d\tau\right] ,\,\dot{x}(\tau,\omega) = f(x(\tau),u(\tau,\upsilon),\omega(\tau))
\]
which, when restricted to explicitly ergodic solutions with stationary state distribution $\rho_x$, becomes
\[
\arg\min\limits_{u(x,\upsilon),\rho_x}\,\int V(x,u) d\rho_x(x),\,\mathcal{P}(x\sim \rho_x)=\{\mathcal{P}(\int_{t}^{t'} f(x,u(x,\upsilon)) dt),\,x\sim \rho_x,\}
\]

See, e.g.~\cite{jenkinson2006ergodic,jenkinson2019ergodic,garibaldi2017ergodic,huang2019ergodic} for references on ergodic optimization and control. 
A standard result is that with a sufficiently regular objective and restriction of the solution to the space of densities with bounded second moments, weak$^*$ compactness holds to ensure (weak$^*$ topology) limits of minimizing sequences $u_k,r_k$. In many cases the solutions are extremal end points of the distribution space. 

In this paper, we consider dynamics that are not long run ergodic. Instead, rather than a single stationary distribution, they are periodic. Thus the probability density is not invariant to a single step of forward propagation, but a number of steps corresponding to its period. 

\subsection{PAR and PAMAR Models}

We follow the presentation and notation given in~\cite{franses2004periodic}. Periodic Time Series models define a stochastic process as a sequential repeated transition defined by an autoregressive model with specific seasons - finite collections of sequential unit time steps organized themselves in a repeating block sequence, e.g. the seasons. Each season will have its own mean and parameters associated with autoregression and moving averages.

Generically, a Periodic Autoregressive Model of order $p$ for stochastic process $Y$ is defined as:
\begin{equation}\label{eq:par}
    Y_t = \mu_{s(t)} + \sum\limits_{i=1}^p \phi_{i,s(t)} Y_{t-i}+\epsilon_t
\end{equation}
where $s(t)$ is the season at time $t$, which can include both the month in the year as well as the hour in the day as a double index. The vector $\mu_{s(t)}$ defines seasonal averages and so this model can represent both auto-regressive influence and the influence of the season's climate. 

In addition to autoregressive influence, we can also incorporate a moving average effect, where again a seasonally dependent parameter can influence the structure of temporal influence in the dynamics. A Periodic Autoregressive Moving Average Model of order $p,q$, with $p$ for the total auto-regressive lag and $q$ the total noise moving average lag, for stochastic process $Y$ is given as:
\begin{equation}\label{eq:pamar}
    Y_t = \mu_{s(t)} + \sum\limits_{i=1}^p \phi_{i,s(t)} Y_{t-i}+\sum\limits_{i=0}^q \theta_{i,s(t)}\epsilon_{t-i}
\end{equation}

In the literature on PAR and PAMAR time series models, various procedures for estimation are available~\cite{franses2004periodic}. In this work we consider the exogenous process defining the relevant periodic inputs as given, and define appropriate schemes to sequentially control such a process in real-time. The target of the control will be to establish the quasi-ergodic dynamics of the system long term with periodicity of similar wavelength strides as the random data. 

Generically, the seasonal index $s(t)$ can have $S$ distinct periodic seasonal patterns that influence the stochastic process. We denote $s(t)=(j_1,j_2,\cdots,j_S)$ where $j_i\in [S_i]$ for $i\in[S]$, that is each periodic process has $S_i$ seasons. We will denote the period for each by $T_i,\, i\in[S]$, with $T\equiv T_S$ the longest period, which will usually be one year.

The use of PAMAR models for applications in the management of energy markets and power grid is standard operating procedure, see e.g.~\cite{conejo2010decision}. This is associated with the intimate relationship of the power grid to weather patterns. For instance the production of energy depends on solar irradiation in the case of solar power and precipitation in the case of hydropower. Likewise, consumption of electricity for turning on lights depends on the ambient level of outdoor light or darkness, and consumption for heating or cooling depends on the temperature outside. Weather patterns in turn depend on planetary geopositioning, and thus exhibits circadian and yearly periodicity. Periodicity in social norms and institutions provide additional drivers, as in the week including distinct consumption patterns on the weekdays and weekend.

\section{Example 1: Hydropower Water Pricing}
Before presenting the generic mathematical form for a control problem for a system exhibiting periodic stochastic dependence in the next Section, we motivate the generic formulation by considering an important application for contemporary management of energy systems. In particular, this Section concerns the problem of Hydropower Water Pricing, which includes managing the opening of dams in a cascade, and thus provision of motive water for power generation, as well as computing the optimal price in the energy markets for this water. 

In general, with the choices of opening dams in a network of rivers and reservoirs, one wishes to avoid spilling while performing discharge when prices are high, as circumstances are such as to encourage matching supply to higher demand. At least since the 1950s \cite{allen1986dynamic}, short-term hydropower scheduling has been approached as an optimal control problem. 

The operation of this optimal control problem, however, is complicated by the uncertainty in regards to important quantities that constrain the hydropower regime or present targets and considerations. There are two fundamental sources of uncertainty when it comes to hydropower scheduling and planning: 1) the weather, especially precipitation, which produces a supply of water into the reservoir, and 2) the demand for electricity, the amount, at varying times, of energy demanded from the grid by connected users, households and businesses.

Clearly there is a significant dependence on the quality of the operation across time periods -- a high water level now presents a greater risk of spill at the next time period, but especially if high precipitation is forecasted. On the other hand, a low water level now with high future demand suggests insufficient supply and thus a welfare-lowering spike in prices to try to balance the market. The operating time period can roll out far into the future, with the interactivity across time periods slowly dissipating, yet always present, as the horizon increases. 

\paragraph{Cost Function}
Consider the simple model for hydropower economic management as defined in~\cite[Chapter 3]{forsund2015hydropower}. 
We adapt it to consider yearly periodicity and the presence of a network of turbines and dams. We aim to maximize social welfare as measured by the consumer surplus, the area under the demand curve,
\[
\sum\limits_{t=1}^T \int_0^{e^H_t} p(z,t) dz
\]
where $e^H_t$ is the total hydropower delivered at time $t$. Note, however, that both the total hydropower delivered and the demand curve can depend randomly, as the delivered hydropower depends on the current water level which depends on past stochastic realizations. When we introduce this explicit randomness, the quantity itself becomes a random variable, and as such it is not uniquely defined. We can impose an expectation or a risk measure, indicating an aversion to negative utility. 

\paragraph{Demand Process Dynamics}

To be precise with the uncertainty, we introduce $(\Omega,\mathcal{T},\mathbb{P})$ to be a probability space. 
Let $\mathcal{F}=\{\mathcal{T}^t\}$ be a 
sequence of sub-sigma algebras of $\mathcal{T}$ such that $\mathcal{T}^t\subset\mathcal{T}^{t+1}$. We let the demand function follow a PAMAR model with seasons $s(t)$ that include both the time of day and time of year,
\begin{equation}\label{eq:demandfun}
    p(z,t) := p_{s(t)}(z)+\sum\limits_{\tau=1}^p \phi^p_{\tau,s(t)} p(z,t-\tau)+\sum\limits_{\tau=0}^q \theta^p_{\tau,s(t)}\zeta_{t-\tau}
\end{equation}
where $\zeta_t$ is a noise term. 

Thus the demand function is itself an auto-regressive moving average stochastic process, a map between price and the quantity of interest to be consumed. 

\paragraph{Decision Variables}

Now let us turn to the controls. The time periods are related in the sense that the amount of water available, which depends on the previous time period, influences the amount to be used or discharged. The total power generated is a function of the controls. As such $e^H_t$ is an intermediate variable and a stochastic state variable. Letting $(\zeta_0,\zeta_1,...,\zeta_t)=\omega_t\in \mathcal{T}^t$ be a random element of the $\sigma$-algebra in $\mathcal{T}^t$, we express the explicit dependence of the production on this variable.

The objective criterion, applying the risk measure to target mitigating the downside risk in revenue, is
defined as:
\begin{equation}\label{eq:obj}
f(e^H_t(\omega_t)):=\sum\limits_{t=1}^T -\mathop{CVaR}\left[-\int_0^{e^H_t(\omega_t)} p(z,t) dz\right]
\end{equation}
Note that we take the standard approximation of ignoring variable costs and only considering revenue as far as the objective function. In a more detailed formulation, we may consider small costs associated with ramping or switching controls. 

\paragraph{State Dynamics}

We turn now to the equations governing the water flow through a hydropower cascade. For the illustrative purposes of this work, we eschew modeling the water dynamics directly, and simply consider that there is a network of turbines indexed by $I$, with some having on-off capacity to let water through and generate power, at some fixed amount, denoted by $u^{(i)}\in\{0,1\}$, and others having both on-off as well as a quantitative control (e.g., RPM of the turbine) $u^{(i)}\in[0,b^u_i]$. Denote the set of feasible controls $\mathbb{U}$, a set of mixed binary variable sets and continuous segments in $\mathbb{R}$. We let the power generated be $e^i_t = a^{(i)} u^{(i)}$ where $a^{(i)}$ is an efficiency rating. The set of equations governing the water flow $r^{(j)}_t$ through a network $j\in J$ is defined by a set of linear equations defining the input of the flow across river networks. 

\begin{equation}\label{eq:waterflowstate}
    r_t(\omega_t) = A r_{t-1}(\omega_{t-1})+Bu_{t-1}+F_t
\end{equation}
where $A$ is a matrix defining the water estuary flow, and $B$ relates the turbine operation to the water in the river sections, and $F_t$ is an exogenous supply (from the original water source and precipitation). While in practice, $F_t$ is uncertain as well, for simplicity we only consider handling the uncertainty in the demand in this work. We do not distinguish between reservoirs and river sections in this formalism, with the arithmetic context sufficient. However, we note that we'd expect the same daily and seasonal periodicity will be exhibited in these quantities. 

Due to capacity and environmental constraints, we could have bound constraints on the water level
\[
\underline{r} \le r_t \le \bar{r}
\]
However, since $r_t$ is a state variable, this is a random quantity, and this needs to be acknowleged in the constraint formulation, which we do so with chance constraints, with a very low probability of violation
\begin{equation}\label{eq:constraints}
\mathbb{P}\left[ \underline{r} \le r_t(\omega_t) \le \bar{r}\right] \ge 1-\alpha_c   
\end{equation}
with $\alpha_c\ll 1$.

We can consider writing the form of an OCP, to be able to solve in an iterative MPC fashion, as follows:
\begin{equation}\label{eq:mpc0}
\begin{array}{rl}
\min\limits_{u\in\mathbb{U}\times \mathcal{T}\times \Xi,r\in \mathbb{R}^{\vert J\vert}\times\mathcal{T}\times \Xi} & \frac{1}{T}\sum\limits_{t=1}^T -\mathop{CVaR}\left[-\int_0^{e^H_t(\omega_t)} p(z,\omega(t)) dz\right]
\\ 
\text{subject to } & r_t(\omega_t) = A r_{t-1}(\omega_{t-1})+Bu_{t-1}(\omega_{t-1})+F_t \\
& e^H_t = \sum_I  e^{(i)}_t =  \sum_I  a^{(i)} u^{(i)}_t(\omega_{t}) \\
& \mathbb{P}\left[ \underline{r} \le r_t(\omega_t) \le \bar{r}\right] \ge 1-\alpha_c   
\end{array}
\end{equation}
where observe that here we make the control depending on the current time realization of the uncertainty, as well as past and predicted stochasticity. This defines the parametrization by $\Xi=\{\omega_t\}_{t=0,\cdots,T}$. We note that standard multistage stochastic programming non-anticipativity can be readily enforced. 

\paragraph{Periodicity Constraint}

However, there is an important missing piece in this formulation: a terminal cost or constraint. Indeed, the purpose of a terminal criterion is to enforce stability of the system. That is, so that the closed loop operation of repeatedly solving a receding horizon OCP results in a reasonably good solution to the infinite (or very long) horizon problem it is meant to approximate. Without any such modification, a vacuous solution to~\eqref{eq:mpc0} can be an objective minimizing by  withdrawing all of the water and leaving the bare minimum at the final time period $T$. However, presumably the operation of the hydropower cascade for time $T+1$ still must occur, and such a myopic policy would not be acceptable. 

There are a number of techniques to address this. Formal terminal cost or constraint criteria can be developed, however this is challenging, especially in the stochastic case. Otherwise some heuristic criteria could approximately suffice. Alternatively a reinforcement learning approach could solve this problem approximately in the case of $T\to \infty$. However, the state space $r$ here is infinite, suggesting such an approach would be intractable. While approximations exist of varying precision~\cite{powell2007approximate}, we choose an alternative approach, noting that we expect the demand to have a cyclical distribution. That is, we expect that the distribution of demand is periodic, repeating yearly.

Recall that to model the noise component, we use an auto-regressive model in~\eqref{eq:demandfun}. We may also incorporate drift into the mean of the noise itself, and define:
\begin{equation}\label{eq:demandnoise}
    \zeta_t = \sum\limits_{\tau=0}^q \theta^q_{\tau,s(t)}\zeta_{t-\tau}+\eta_t,\,\eta_t\in \mathcal{D}(\bar{\eta}_t),\,\bar{\eta}_t = \bar{\eta}_{t+T},\,\forall t
\end{equation}
where $\mathcal{D}$ is a sub-Gaussian process with variance $\sigma$ and $\bar{\eta}_t$ is an exogenous mean demand function that is cyclical with period of one year, which we denote now $T^Y$.

Now consider taking $T\to \infty$ in~\eqref{eq:mpc0}, yielding an infinite sum and set of constraint criteria. However, by considering the asymptotic scenario, initial conditions and transient behavior can be seen to have negligible effect on the long run cost if the system converges, in distribution, to a process of periodic ergodicity, that is, if $\omega(t)$ is such that $\omega(t)=\omega(t+T^Y)$ then,
\[
\lim\limits_{T\to \infty} \frac{1}{T}\sum\limits_{t=1}^T -\mathop{CVaR}\left[-\int_0^{e^H_t(\omega_t)} p(z,\omega(t)) dz\right] \approx \frac{1}{T^Y}\sum\limits_{t=1}^{T^Y} -\mathop{CVaR}\left[-\int_0^{e^H_t(\omega_t)} p(z,\omega(t)) dz\right]
\]
and so we can forego the decicate and fraught choice of having to decide on a horizon time period length, recognizing that we can just focus on optimizing the long run average by considering one period of the cyclical stochastic process, as driven by one period of a cycle in the distribution of precipitation and demand, for which the natural choice is one year. 

This suggests the notion associated to stability in consider \emph{periodic ergodicity}. In particular, the novel technical tool in this work is an explicit enforcement of probabilistic cyclicity. Formally, with $\mathcal{P}(\cdot)$ denoting the probability law of a random variable:
\begin{equation}\label{eq:ergcon}
r_t(\omega_t) = \sum\limits_{i=1}^S r_t^{(i)}(\omega_t)
    ,\,\,\mathcal{P}(r_t^{(i)}(\omega_t)) = \mathcal{P}(r_t^{(i)}(\omega_{t+T_i})),\,i\in[S]
\end{equation}
Thus, we decompose the state into a sum of periodic stochastic processes, with periods coinciding with those of the exogenous noise seasons.

We now define the full-space ergodic control problem with multi-periodicity enforcement:
\begin{equation}\label{eq:probdesign}
\begin{array}{rl}
\min\limits_{u_t\in\mathbb{U}\times \Xi,r_t\in \mathbb{R}^{\vert J\vert}\times \Xi,t\in\mathcal{T}} & \sum\limits_{t=1}^T -\mathop{CVaR}\left[-\int_0^{e^H_t(\omega_t)} p(z,\omega(t)) dz\right]
\\ 
\text{subject to } & r_t(\omega_t) = A r_{t-1}(\omega_{t-1})+Bu_{t-1}(\omega_{t-1})+F_t \\
& e^H_t = \sum_I  e^{(i)}_t =  \sum_I  a^{(i)} u^{(i)}_t(\omega_{t}) \\
& \mathbb{P}\left[ \underline{r} \le r_t(\omega_t) \le \bar{r}\right] \ge 1-\alpha_c \\
&     r_t(\omega_t) = \sum\limits_{i=1}^S r_t^{(i)}(\omega_t)
    \\ & \mathcal{P}(r_t^{(i)}(\omega_t)) = \mathcal{P}(r_i(\omega_{t+T_i})),
\forall t\in\mathcal{T},\, i\in[S]
\end{array}    
\end{equation}

Note that this problem is still in its fully probabilistic form, with the states and controls dependent on the possibly infinite space the random variables influencing the problem depend on. This permits a more careful and principled choice of discretization further down the implementation pipeline, as well as more flexible algorithms with repeated stochastic discretizations. 

However, we need to understand the appropriate definition and modeling of $u(\omega_0)$, that is, the (randomized) choice at the beginning of the year. In particular, without treating the past explicitly, this would imply the first control is entirely uninformed, that is beyond the observed value at the initial state, i.e., $\omega_0$ is associated with realizations of, only, $r_0$. Of course, unless we are working with a newly operational station, this would be needlessly discarding available information. On the other hand if we explicitly include dependence on past years, the associated sigma algebra cylinder still is monotonically expanding, thus no longer creating a natural symmetry across years.

This can be entirely mitigated by restricting ourselves to \emph{deterministic} policies. This would mean the functional dependence of the decision variables will be $u(t,r_t)$. We will consider that the PAMAR model of the distribution is taken as a given, and we first consider the problem without stochastic transience. We consider relaxing this in the next subsection, with the ultimate operating procedure solving a precise offline solution that will serve as a set point for real-time solvers. 

The final form of the problem, with the decision variable properly defined, and the statistics of $\{\omega_t\}$ are taken to be the periodically ergodic distribution of~\eqref{eq:pamar}:
\begin{equation}\label{eq:ergstatdet}
\begin{array}{rl}
\min\limits_{u(t,x)\in\mathbb{U},r_t\in \mathbb{R}^{\vert J\vert}\times \Xi,\, t\in \mathcal{T}} & \sum\limits_{t=1}^T -\mathop{CVaR}\left[-\int_0^{e^H_t(r_{t-1}(\omega_{t-1}))} p(z,\omega(t)) dz\right]
\\ 
\text{subject to } & r_t(\omega_t) = A r_{t-1}(\omega_{t-1})+Bu_{t-1}+F_t \\
& e^H_t = \sum_I  e^{(i)}_t =  \sum_I  a^{(i)} u^{(i)}_t \\
& \mathbb{P}\left[ \underline{r} \le r_t(\omega_t) \le \bar{r}\right] \ge 1-\alpha_c \\
&     r_t(\omega_t) = \sum\limits_{i=1}^S r_t^{(i)}(\omega_t)
    \\ & \mathcal{P}(r_t^{(i)}(\omega_t)) = \mathcal{P}(r_t^{(i)}(\omega_{t+T_i})),
\forall t\in\mathcal{T},\, i\in[S]
\end{array}    
\end{equation}
where $u_t$ can be a policy that is a parametrized function $u_t=f(r_t;t)$ that depends on the realization of the state at time $t$, $\hat{r}_t$. The term $\Pi$ is a projection of the cross seasonal dynamics' moving averages. That is, the density on Monday one week should be the same as the density the next week, modulo steady increase as the weather gets colder, which is defined by the yearly periodic terms.






\paragraph{Transient Control Problem}
It can be seen that~\eqref{eq:ergstatdet} is a problem that depends only on the forecasted probabilities and can be solved offline. Thus, independent of any real time of implementation, the solution of the mathematical form of this ergodic control problem can be computed and catalogued. This permits a fine-grained discretization in time and stochastic dimension, with the solution obtained during a long time operation on an HPC cluster. 

However, in real time, of course, the actual noise of the process, that is the weather and influence of precipitation and demand, obtains a realization. Thus a more refined solution can be obtained in real time that is with respect to the actual recent past of stochastic realizations and the history-depending forecast in the near future. Of course, with real time implementation, the same problem size dimension cannot be handled, rather, for one the discretization must be coarser. However, with solving a transient problem, the entire year cycle time horizon does not need to be considered. Instead, a far shorter time length $\check{T}\ll T$ can be applied, and as a terminal state, the solution to the offline problem~\eqref{eq:ergstatdet}. 

Consider having computed the offline solution to~\eqref{eq:ergstatdet} above as $\{u^*(t,x),r^*_t\}$. One is to perform a real time control at a particular time $\bar{t}$. That is, with $\omega_t$ now realized in a given year up to $\bar{t}$, denoted $\{\hat{\omega}_{t}\}$ with $t\le\bar{t}$ and $\check{\omega}_t$ the forecasted future noise computed by propagating the statistical time series model~\eqref{eq:pamar} forward until $\check{T}$, we solve,

\begin{equation}\label{eq:ergstattran}
\begin{array}{rl}
\min\limits_{u(t,x)\in\mathbb{U},r_t\in \mathbb{R}^{\vert J\vert}\times \Xi,\, t\in \{\bar{t},\cdots,\check{T}\}} & \sum\limits_{t=\bar{t}}^{\check{T}} -\mathop{CVaR}\left[-\int_0^{e^H_t(r_{t-1}(\check{\omega}_{t-1}))} p(z,\check{\omega}(t)) dz\right]
\\ 
\text{subject to } & r_t(\check{\omega}_t) = A r_{t-1}(\check{\omega}_{t-1})+Bu_{t-1}+F_t \\
& e^H_t = \sum_I  e^{(i)}_t =  \sum_I  a^{(i)} u^{(i)}_t \\
& \mathbb{P}\left[ \underline{r} \le r_t(\check{\omega}_t) \le \bar{r}\right] \ge 1-\alpha_c \\
&     \mathcal{P}(r(\omega_{\bar{t}})) = \delta(\hat{r}_{\bar{t}})\\
& \mathcal{P}(r(\check{\omega}_{\check{T}})) = \mathcal{P}(r^*(\omega_{\check{T}})),
\end{array}    
\end{equation}
then implement the result $\check{u}(\bar{t})$ and then proceed to the next time step, formulating the transient problem again. Observe that the final state is directed to the offline computed solution. 



\section{Generic Problem Statement}
Here we define a generic form of the multi-periodic process control problem. We begin with the full period (typically a year) offline problem and then the transient. We generically define:
\begin{enumerate}
    \item The decision variable in some appropriate space $\mathbb{U}$ over time set $\mathcal{T}=\{1,\cdots,T\}$ and current state $x$, denoted $u(t,x)$
    \item The state variable defined in an appropriate space $\mathcal{X}$ over the same time set $\mathcal{T}$ and stochasticity $\Xi=\{\omega_1,\cdots,\omega_T\}$.
    \item Stochastic process $\omega_t$ modeled by a multi-periodic PAMAR process defined in~\eqref{eq:pamar} with periods $\{T_1,\cdots,T_S\}$.
    \item We use $\mathcal{R}$ to generically denote a statistical aggregation for the cost function, which can be
    \item We let $f$ define the overall time dynamics of the state variable given the previous state, control, and random noise
    \item We let $V$ denote the form of the objective function itself.
    \item The function $h$ will define an inequality constraint that we will write as a chance constraint. Note, however, that more computationally straightforward approaches, such as a barrier function ensuring almost sure satisfaction, may be used instead.
    \item The functional $\mathcal{R}$ will represent a generic statistical aggregation, which can correspond to an expectation or a risk measure. 
\end{enumerate}

\begin{equation}\label{eq:ergstatdetgeneral}
\begin{array}{rl}
\min\limits_{u(t,x)\in\mathbb{U},x_t\in \mathbb{R}^{\vert J\vert}\times \Xi,\, t\in \mathcal{T}} & \sum\limits_{t=1}^T \mathcal{R}\left[V(x_t(\Xi),\omega_t,u_t(x_{t-1}(\Xi)))\right]
\\ 
\text{subject to } & x_t(\Xi) = f(x_{t-1}(\Xi),\omega_t,u_t(x_{t-1})) \\

& \mathbb{P}\left[ \underline{h} \le h(x_t(\Xi)) \le \bar{h}\right] \ge 1-\alpha_c \\
&     x_t(\omega_t) = \sum\limits_{i=1}^S x_t^{(i)}(\omega_t)
    \\ & \mathcal{P}(x_t^{(i)}(\omega_t)) = \mathcal{P}(x_t^{(i)}(\omega_{t+T_i})),
\forall t\in\mathcal{T},\, i\in[S]

\end{array}    
\end{equation}

Now we have defined the Optimal Control problem corresponding the multi-periodic ergodic optimization. 

There are a number of questions that this introduces: as far as the theoretical properties, in particular the existence, uniqueness and regularity of solutions.

\subsection{Properties}
The problem defined above is general, and without specific assumptions regarding the dynamics $f$ and the constraint function $h$, one cannot guarantee that a solution $\{u(t,x)\}$ exists such that there exists a stochastic process defining the state $\{x_t(\Xi\}\}$ that satisfies all of the constraints. It is thus important to perform accurate identification of these dynamics and ensure that the range of $f$ under different feasible control actions is sufficiently regular as well as explorative of the state space. 

For the problem applications discussed in this paper, it is naturally expected that the state space and controls will be permissible from a bounded set. This can be ensured with a bounded control constraint set $\mathbb{U}$ and a contractive or Lyapunov-associated objective $f$~\cite{meyn2012markov}. With a lower semicontinuous $V$ and convex $\mathcal{R}$, this should be sufficient for establishing the existence of a solution by standard arguments, assuming that the chance constraint isn't onerous.

\subsection{Practical Considerations for Solution Computation}

as otherwise $u(x)$ is an infinite-dimensional object. 
Most commonly, Sample Average Approximation (SAA) sampling of scenario trees is the standard approach to these problems. The text~\cite{conejo2010decision} presents a comprehensive guide to generating scenarios that follow various Autoregressive structure, with PAMAR models taking a prominent role. A number of practical and technical considerations are discussed for the efficient and useful generation of scenario samples.

SAA methods introduce significant scaling problems as due to their exponential complexity. In this case, the advantage of the formulation above is that we can specifically define a basis, using a spectral, hyperspectral, or wavelet expansion, to define,
\[
x_t(\omega_t)\approx \sum\limits_{i=1}^S\sum\limits_{m=1}^{M_i} \alpha_{t,i,m}\psi_{i,m}(t,x)
\]
with $\psi_{i,m}(t,x)$ are basis vectors with period, in $t$, of $T_i$. A natural approach would be to define,
\[
\psi_{i,m}(t,x) = e^{-ir_m t} \phi_m(x)
\]
where $e^{-r_mt}$ is a Fourier expansion basis element and $\phi_m(x)$ a basis element of a polynomial chaos expansion over $x$. 

By using the numerical stochastic approximation that explicitly enforces the periodic ergodicity constraints, one can provide an immediate structural ansatz satisfying these constraints and focus on optimizing the criteria with the dynamics and chance constraint. Furthermore, we can discretize the rest of the problem by evaluating the functions at the basis vectors and optimizing over coefficients. 

In the case of linear dynamics $f$, standard polynomial chaos expansions lead to straightforward sets of linear systems of equations to solve, as far the constraint discretization. Otherwise, generalized polynomial chaos is used to facilitate computation through the pushforward of the nonlinear map $f$, see, e.g.~\cite{janya2017framework}.
\section{Additional Examples}
\subsection{Example 2: Virtual Power Plant Operation}

We write the following Robust Virtual Power Plant Electricity Market problem from~\cite[Chapter 4, Section 4.3.3]{baringo2020virtual}

\begin{equation}\label{eq:ergstatdetgeneral}
\begin{array}{rl}
\min\limits_{\Phi_t\in H,\, t\in \mathcal{T}} & \sum\limits_{t=1}^T \mathbf{CVaR}\left[\pi(\Phi_t(\omega_t),\omega_t)\right]
\\ 
\text{subject to } & \Phi^x_{t}(\omega_t) = f(\Phi^x_{t-1}(\omega_{t-1}),\omega_t,\Phi^u_t(\Phi^x_{t-1})) \\

& \underline{h} \le \Phi^x_{t}(\omega_t) \le \bar{h} \text{  almost surely}\\
& \Phi^u_t(\Phi^x_t)\in\mathbb{U} \\
&     \Phi^x_t(\omega_t) = \sum\limits_{i=1}^S \Phi_t^{x,(i)}(\omega_t)
    \\ & \mathcal{P}(\Phi_t^{x,(i)}(\omega_t)) = \mathcal{P}(\Phi_t^{x,(i)}(\omega_{t+T_i})),
\forall t\in\mathcal{T},\, i\in[S]
\end{array}    
\end{equation}
where we represent $\Phi_t=(\Phi_t^x,\Phi^u_t)$ as the state and control variables in a compact form. 

The control vectors in this Example include, for each time period:
\begin{enumerate}
    \item Energy drawn down from battery unit
    \item Power generation from conventional plant
    \item Available capacity for renewable generation
    \item Power traded in the Day Ahead Energy Market
\end{enumerate}
and the state variables are:
\begin{enumerate}
    \item Energy stored in battery unit
    \item Stochastic renewable power generation
    \item Power flow through each line
    \item Power traded in the real time energy market
    \item Charging up of the battery unit
\end{enumerate}
The cost function $\pi$ is the negative profit, the dynamics $f$ define the power flow and market prices, and the constraints enforce market clearing, i.e., a matching of energy production to demand. 

In the text~\cite{baringo2020virtual}, three time periods are considered, with four scenarios generated at each stage. The definition of the constraints for all of the scenarios is demonstrated, and a solution computed and presented.

Observing that a more realistic scenario would require a finer stochastic discretization, while at the same time, one should note that in many cases, conventional power plants, or pumped hydro energy, can be provided in case of a shortfall in supply, however it must be ramped up beforehand. By considering more long run scenarios, the VPP can incorporate proper foresight in its operation.

The constraints in~\cite{baringo2020virtual} are linear, making computation of the solution coefficients in a discretization straightforward. 
\subsection{Example 3: Agriculture Management}
Agricultural management, that is the decisions of planting and cutting and quantity of nutrients and irrigation to supply, presents a challenging problem of long run horizon control, with phenomena defined by cyclical weather patterns. See~\cite{ding2018model,hertzler1991dynamic} for some expositions of model prediction control and stochastic control for agriculture.

\begin{equation}\label{eq:ergstatdetgeneral}
\begin{array}{rl}
\min\limits_{\Phi_t\in H,\, t\in \mathcal{T}} & \sum\limits_{t=1}^T \mathbf{CVaR}\left[\pi(c_t(\omega_t),a_t(c_t),\omega_t)\right]
\\ 
\text{subject to } & c_{t}(\omega_t) = f(c_{t-1}(\omega_{t-1}),\omega_t,a_t,c_{t-1})) \\

& \underline{h} \le c_{t}(\omega_t) \le \bar{h} \text{  almost surely}\\
& a_t(\Phi^x_t)\in\mathbb{A}\subset\mathbb{R}^{a_c}\times \{0,1\}^{a_b} \\
&     c_t(\omega_t) = \sum\limits_{i=1}^S c_t^{(i)}(\omega_t)+c^{(0)}_t
    \\ & \mathcal{P}(c_t^{(i)}(\omega_t)) = \mathcal{P}(c_t^{x,(i)}(\omega_{t+T_i})),
\forall t\in\mathcal{T},\, i\in[S]
\end{array}    
\end{equation}
Here $a_t$ consists of $a_c$ continuous decisions (water to pour, mass to cut), and $a_d$ discrete decisions (planting a seed, turning on an irrigation system). The vector of crop biomass $c_t$ experiences both a self driven natural growth $c^{(0)}_t$ together with the periodic stochastic terms that slow or speed up growth due to additional sunlight or precipitation, etc. 

The vector $\pi$ represents the profit obtained by the farmer from selling the production. The prices as well as the production are driven by cyliccal daily and yearly patterns. We can establish a production system that presents a robust operation of stabilizing returns across the year, through appropriate adaptive management decisions.

\section{Conclusion and Ongoing Work}

This paper presented a formulation for controlling a system governed by periodic autoregressive state dynamics. These are a common occurrence with production systems that depend on the weather and climate, and the structure of the noise presents an opportunity to obviate the usual scenario generation complexity explosion by encoding explicitly periodic state dynamics in probability. Natural future work is testing the method numerically as well as comparing it to standard alternatives, as well as seeking to obtain specific optimization approximation and stability guarantees.





\paragraph*{Acknowledgements}
This work received funding from the National Centre for Energy II (TN02000025).

\bibliographystyle{plain}
\bibliography{refs}

\clearpage 
\appendix

\end{document}